\documentclass[12pt]{article}

\usepackage{amsfonts}

\newtheorem{thm}{Theorem}[section]

\newtheorem{cor}[thm]{Corollary}

\newtheorem{prop}[thm]{Proposition}

\newenvironment{remark}{\par\medskip\noindent{\bf Remark.\ }}{\par\smallskip}
\newcommand{\proof
}{\par\medskip\noindent {\bf Proof.\ \ }}

\newcommand{\be}{\begin{equation}}
\newcommand{\ee}{\end{equation}}
\newcommand{\openbox}{\leavevmode
  \hbox to8pt{\hfil\vrule\vbox to6pt{\hrule width6pt\vfil\hrule}\vrule}}

\newcommand{\qed}{\hbox to5pt{ } \hfill \openbox\bigskip\medskip}

\newcommand{\cP}{\mbox{$\cal P$}}

\newcommand{\Z}{\mathbb Z}
\newcommand{\Q}{\mathbb Q}
\newcommand{\R}{\mathbb R}

\title{The surface of a lattice polytope}
\author{G\'abor Heged\"{u}s
\\{\normalsize Johann Radon Institute for Computational and Applied Mathematics}
}

\begin{document}

\footnotetext{
{\bf Keywords.} Lattice polytopes, Ehrhart polynomial 

{\bf 2010 Mathematics Subject Classification.} 52B20,52C07, 11H06 }

\maketitle

\begin{abstract}
My main results are simple formulas for the surface area of $d$--dimensional lattice polytopes using Ehrhart theory. 
\end{abstract}
\medskip

\section{Introduction}
\noindent

Throughout the paper a {\em lattice polytope} $P\subseteq \R^d$ is a polytope whose vertices have integral coordinates.
 
Let $S\subseteq \R^d$ be a subset of the Euclidean space $\R^d$. Let $G(S)$ denote the {\em lattice point enumerator} of the set $S$, the number of lattice (integral) points in $S$, i.e., $G(S)=|(S\cap \Z^d)|$.  

Let $\cP$ denote an arbitrary $d$--dimensional lattice polytope. In the following we denote by 
$$
\nu \cP:=\{n\underline{x}:~ \underline{x}\in P\}
$$
the dilatate of $\cP$ by the integer factor $\nu \geq 0$. 

In 1962 E. Ehrhart proved (see e.g. \cite[Chapter 3, Chapter 5]{BR}, \cite{E}) the following Theorem:

\begin{thm} \label{Ehr}
Let $\cP$ be a convex $d$--dimensional lattice polytope in the Euclidean space $\R^d$. 
Then there exists a unique polynomial (the {\em Ehrhart polynomial})
$$
E_{\cP}(x):=\sum_{i=0}^d e_i(\cP)x^i\in \Q[x],
$$
which has the following properties:\\
(1) For all integers $\nu\geq 0$,
$$
E_{\cP}(\nu)=|(\nu \cP)\cap \Z^d|.
$$
(2) The leading coefficient $e_d(\cP)$ of $E_{\cP}(x)$ is $\mbox{vol}(\cP)$, the volume of $\cP$. \\
(3) If $\mbox{int}(\cP)$ denotes the interior of $\cP$, then the {\em reciprocity law} states that for all integers $\nu >0$,
\begin{equation}
E_{\cP}(-\nu)=(-1)^d |(\nu\cdot {\rm int}(\cP))\cap \Z^d|.
\end{equation}
(4) The second leading coefficient $e_{d-1}(\cP)$ of $E_{\cP}(x)$ is the half of the lattice surface area of $\cP$:
$$
e_{d-1}(\cP)=\frac{1}{2}\sum_{F\ facet\ of\ \cP} \frac{{\rm vol}_{d-1}(F)}{{\rm det}({\rm aff}\ F\cap \Z^d)}.
$$
Here $\mbox{vol}_{d-1}(\cdot)$ denotes the $(d-1)$--dimensional volume and $\mbox{det}(\mbox{aff }F\cap\Z^d)$ denotes the determinant of the $(d-1)$--dimensional sublattice contained in the affine hull of $F$. \\
(5) The constant coefficient $e_{0}(\cP)$ of $E_{\cP}(x)$ is $1$.
\end{thm}

Let $\cP$ be a convex $d$--dimensional lattice polytope, which contains the origin of the lattice in its interior.
We say that $P$ is {\em reflexive} if the dual polytope $P^*$ is a lattice polytope, where the {\em dual polytope} of $P$ is defined as
$$
P^*:=\{y\in {\R}^d:~ \langle x,y \rangle \geq -1 {\rm\ for\ all\ } x\in P\}.
$$

A. M. Kasprzyk proved in \cite[Proposition 3.9.2]{K} the following equivalent characterization of reflexive Fano polytopes:
\begin{prop} \label{volume}
Let $\cP$ be a $d$--dimensional Fano polytope. Then $\cP$ is reflexive iff 
\begin{equation}
{\rm vol}(\cP)=\frac{{\rm surf} \cP}{d}.
\end{equation}
\end{prop}

In 1899 G. A. Pick published his famous formula in \cite{P}. Using this formula we can compute easily the area of a lattice polygon.  
Pick showed that the following expression gives the area of a simple lattice polygon $Q$:
$$
\mbox{Area}(Q)=I+\frac{B}{2}-1,
$$
where $B$ is the number of lattice points on the boundary of $Q$ and $I$ is the number of lattice points in the interior of $Q$. 

This formula can be derived easily from Ehrhart Theorem \ref{Ehr} (see e.g. \cite[Chapter 4]{C}).

My main results are similar simple formulas for the surface area of $3$--dimensional and $4$--dimensional lattice polytopes using Ehrhart theory.  

\section{The main results}

Let $\cP$ be a convex $d$--dimensional lattice polytope in $\R^d$. Denote by $\mbox{surf}(\cP)$ the lattice surface area of $\cP$:
$$
\mbox{surf}(\cP):=\sum_{F\ facet\ of \cP} \frac{\mbox{vol}_{d-1}(F)}{\mbox{det}(\mbox{aff }\ F\cap \Z^d)}.
$$

Let $i(\cP)$ and $b(\cP)$ denote the numbers $|\mbox{int}(\cP)\cap \Z^d|$ and $|\partial(\cP)\cap \Z^d|$, respectively.

Here $\partial(\cP)$ denotes the boundary of the polytope $\cP$.

\begin{thm} \label{Ehr5} 
Let $\cP$ be a convex $d$--dimensional lattice polytope in $\R^d$.

Suppose that $d$ is an odd number. Let $t:=\frac{d-1}{2}$. Then define the matrix
\[
\mathbf{A}(\cP,d):=\left( 
\begin{array} {cccc}
b(\cP)-2 & 1^{d-3} & \cdots & 1^2 \\ 
b(2\cP)-2 & 2^{d-3} & \cdots & 2^2 \\
\vdots & \vdots & \ddots & \vdots \\
b(t\cP)-2 & t^{d-3} & \vdots & t^2 \\
\end{array} \right) 
\]
and 
\[
\mathbf{D}(\cP,d):=\left( 
\begin{array} {cccc}
1^{d-1} & 1^{d-3} & \cdots & 1^2\\
\vdots & \vdots & \ddots & \vdots \\
t^{d-1} & t^{d-3} & \cdots & t^2\\
\end{array} \right) 
\]

Then
\begin{equation} \label{formprtl}
{\rm surf}(\cP)=\frac{{\rm det}(\mathbf{A}(\cP,d))}{{\rm det}(\mathbf{D}(\cP,d))}
\end{equation}

Suppose that $d$ is an even number.

Let $t:=\frac{d}{2}$. Then define the matrix
\[
\mathbf{B}(\cP,d):=\left( 
\begin{array} {cccc}
b(\cP) & 1^{d-3} & \cdots & 1 \\ 
b(2\cP) & 2^{d-3} & \cdots & 2 \\
\vdots & \vdots & \ddots & \vdots \\
b(t\cP) & t^{d-3} & \vdots & t \\
\end{array} \right) 
\]
and 
\[
\mathbf{D}(\cP,d):=\left( 
\begin{array} {cccc}
1^{d-1} & 1^{d-3} & \cdots & 1 \\
\vdots & \vdots & \ddots & \vdots \\
t^{d-1} & t^{d-3} & \cdots & t \\
\end{array} \right) 
\]

Then
\begin{equation} \label{formps}
{\rm surf}(\cP)=\frac{{\rm det}(\mathbf{B}(\cP,d))}{{\rm det}(\mathbf{D}(\cP,d))}
\end{equation}
\end{thm}

\proof
Let 
$$
E_{\cP}(x):=\sum_{i=0}^d e_i(\cP)x^i\in \Q[x],
$$
denote the Ehrhart polynomial of the polytope $\cP$.

First suppose that $d$ is an odd number. Let $0\leq k\leq \frac{d-1}{2}$. Then
using Theorem \ref{Ehr} 
\begin{equation} \label{Ehrpol}
i(k\cP)+b(k\cP)=L_{\cP}(k)=e_d(\cP)k^d+e_{d-1}(\cP)k^{d-1}+\ldots +1
\end{equation}
and 
\begin{equation} \label{Ehrpol2}
-i(k\cP)=-L_{\cP}(-k)=-e_d(\cP)k^d+e_{d-1}(\cP)k^{d-1}-\ldots +1
\end{equation}
Suming (\ref{Ehrpol}) and (\ref{Ehrpol2}) we get that
$$
b(k\cP)=2e_{d-1}(\cP)k^{d-1}+2e_{d-3}(\cP)k^{d-1}+\ldots +2
$$
i.e.,
$$
b(k\cP)-2=2e_{d-1}(\cP)k^{d-1}+2e_{d-3}(\cP)k^{d-3}+\ldots 2e_{2}(\cP)k^{2}
$$
for each $0\leq k\leq \frac{d-1}{2}$. Solving this linear equation system using Cramer's rule we get that
$$
e_{d-1}(\cP)=\frac{{\rm det}(\mathbf{A}(\cP,d))}{2{\rm det}(\mathbf{D}(\cP,d))}.
$$
But using Theorem \ref{Ehr} (4) we get that
$$
{\rm surf}(\cP)=\frac{e_{d-1}(\cP)}{2},
$$
and we get our result for odd $d$.

Suppose that $d$ is an even number. 
Let $t:=\frac{d}{2}$. Let $0\leq k\leq \frac{d}{2}$. 
Using Theorem \ref{Ehr} 
\begin{equation} \label{Ehrpol3}
i(k\cP)+b(k\cP)=L_{\cP}(k)=e_d(\cP)k^d+e_{d-1}(\cP)k^{d-1}+\ldots +1
\end{equation}
and 
\begin{equation} \label{Ehrpol4}
i(k\cP)=L_{\cP}(-k)=e_d(\cP)k^d-e_{d-1}(\cP)k^{d-1}+\ldots +1
\end{equation}
Substracting (\ref{Ehrpol4}) from (\ref{Ehrpol3}) we get that 
\begin{equation} \label{Ehrpol5}
b(k\cP)=2e_{d-1}(\cP)k^{d-1}+2e_{d-3}(\cP)k^{d-1}+\ldots +2e_1(\cP)k
\end{equation}
for each $0\leq k\leq \frac{d}{2}$. We can again solve this linear equation system using Cramer's rule, hence
\begin{equation} \label{coeff1}
e_{d-1}(\cP)=\frac{{\rm det}(\mathbf{B}(\cP,d))}{2{\rm det}(\mathbf{D}(\cP,d))}.
\end{equation}
Theorem \ref{Ehr} (4) implies that
\begin{equation} \label{coeff2}
{\rm surf}(\cP)=\frac{e_{d-1}(\cP)}{2},
\end{equation}
and we get our result from (\ref{coeff1}) and (\ref{coeff2}).

{\bf Examples.} \\

If $d=3$, then ${\rm surf}(\cP)=b(\cP)-2$. \\

In \cite[Proposition 10.3.2]{K} A. M. Kasprzyk proved this formula from Pick's Theorem.

If $d=4$, then 
$$
{\rm surf}(\cP)=\frac{b(2\cP)-2b(\cP)}{6}.
$$ 
If $d=5$, then 
$$
{\rm surf}(\cP)=\frac{b(2\cP)-4b(\cP)-6}{12}
$$ 

\begin{remark}
A. M. Kasprzyk called my attention for the following consequence of Theorem \ref{Ehr5}.
\end{remark}

\begin{cor} \label{Fano3}
Let $\cP$ be a convex $d$--dimensional Fano lattice polytope in $\R^d$.

Suppose that $d$ is an odd number. Let $t:=\frac{d-1}{2}$. Then define the matrix
\[
\mathbf{A}(\cP,d):=\left( 
\begin{array} {cccc}
b(\cP)-2 & 1^{d-3} & \cdots & 1^2 \\ 
b(2\cP)-2 & 2^{d-3} & \cdots & 2^2 \\
\vdots & \vdots & \ddots & \vdots \\
b(t\cP)-2 & t^{d-3} & \vdots & t^2 \\
\end{array} \right) 
\]
and 
\[
\mathbf{D}(\cP,d):=\left( 
\begin{array} {cccc}
1^{d-1} & 1^{d-3} & \cdots & 1^2\\
\vdots & \vdots & \ddots & \vdots \\
t^{d-1} & t^{d-3} & \cdots & t^2\\
\end{array} \right) 
\]

Then
\begin{equation} \label{formprtl3}
{\rm vol}(\cP)=\frac{{\rm det}(\mathbf{A}(\cP,d))}{d\cdot{\rm det}(\mathbf{D}(\cP,d))}
\end{equation}

Suppose that $d$ is an even number.

Let $t:=\frac{d}{2}$. Then define the matrix
\[
\mathbf{B}(\cP,d):=\left( 
\begin{array} {cccc}
b(\cP) & 1^{d-3} & \cdots & 1 \\ 
b(2\cP) & 2^{d-3} & \cdots & 2 \\
\vdots & \vdots & \ddots & \vdots \\
b(t\cP) & t^{d-3} & \vdots & t \\
\end{array} \right) 
\]
and 
\[
\mathbf{D}(\cP,d):=\left( 
\begin{array} {cccc}
1^{d-1} & 1^{d-3} & \cdots & 1 \\
\vdots & \vdots & \ddots & \vdots \\
t^{d-1} & t^{d-3} & \cdots & t \\
\end{array} \right) 
\]

Then $\cP$ is a reflexive polytope iff 
\begin{equation} \label{formps3}
{\rm vol}(\cP)=\frac{{\rm det}(\mathbf{B}(\cP,d))}{d\cdot{\rm det}(\mathbf{D}(\cP,d))}.
\end{equation}
\end{cor}

\qed
\proof

Corollary \ref{Fano3} is the obvious consequence of Theorem \ref{Ehr5} and Proposition \ref{volume}. \qed

{\bf Acknowledgements.} 
I am indebted to Josef Schicho and Alexander M. Kasprzyk for their useful remarks.

\end{document}